\documentclass[12pt,letterpaper]{amsart}
\usepackage{amssymb,latexsym, amsmath, amsxtra}
\usepackage[dvips]{graphics}

\theoremstyle{plain}
\newtheorem{theorem}{Theorem}[section]
\newtheorem{corollary}[theorem]{Corollary}
\newtheorem{conjecture}[theorem]{Conjecture}
\newtheorem{lemma}[theorem]{Lemma}
\newtheorem{proposition}[theorem]{Proposition}
\theoremstyle{definition}

\theoremstyle{remark}

\numberwithin{equation}{subsection}
\numberwithin{theorem}{subsection}

\newcommand{\ve}{\varepsilon}
\newcommand{\R}{\mathbb R}
\newcommand{\Z}{\mathbb Z}
\newcommand{\Q}{\mathbb Q}

\def\({\left(}
\def\){\right)}

\newcommand{\ontop}[2]{\genfrac{}{}{0pt}{}{#1}{#2}}

\begin{document}
\title[Differentiation Evens Out Zero Spacings]
{Differentiation Evens Out Zero Spacings }
\author{David W. Farmer and Robert C. Rhoades}

\thanks{Research of the first author supported by the
American Institute of Mathematics
and the NSF}

\date{\today}
\thispagestyle{empty}
\vspace{.5cm}
\begin{abstract}
If $f$ is a polynomial with all of its roots on the
real line, then the roots of the derivative $f'$ are more evenly
spaced than the roots of $f$.
The same holds for
a real entire function of order~1 with all its zeros on a line.
In particular, we show that if $f$ is entire of order~1 and has sufficient
regularity in its zero spacing, then under repeated differentiation
the function approaches, after normalization, the cosine function.
We also study polynomials with all their zeros on a circle, 
and we find a close analogy between the two situations.
This sheds light on the spacing between zeros of the Riemann zeta-function 
and its connection to random matrix polynomials.
\end{abstract}

\address{
{\parskip 0pt
American Institute of Mathematics\endgraf
farmer@aimath.org\endgraf
\null
Bucknell University\endgraf
rrhoades@bucknell.edu\endgraf
}
  }

\maketitle

\section{Introduction}

If $f$ is a polynomial with all of its zeros on the real line, then
all of the zeros of the derivative $f'$ also lie on the real line.
The same holds for entire functions of order~1 which are real on the real line.
But while the zeros of $f$ in an interval could in principle be any finite subset
of the interval, 
the condition that $f$ have only real zeros will impose restrictions
on the zeros of~$f'$.  
The main idea of this paper is that
if the zeros of $f$ lie on a line then the zeros of $f'$ are 
more evenly spaced than the zeros of~$f$.
Repeated differentiation leads the zeros to become more and more 
evenly spaced, and under appropriate conditions the zeros will
approach equal spacing.  See Section~\ref{sec:longterm} for precise statements.

There is a long history to the behavior of the zeros of entire functions
under repeated differentiation.  Polya~\cite{p1} made several
interesting conjectures for the case of real entire functions, some
of which have only been solved recently~\cite{CCS1, CCS2,kikim,K1,kim,ss}.
Many of those questions concern the reality of zeros of a real
entire function, or the motion of zeros toward (or away from)
the real line when the function is differentiated.  In this
paper we take a complementary view and consider the motion along
the real line when the function is differentiated.

\subsection{Zeros of the derivative}
\begin{figure}[htp] 
\begin{center}
\scalebox{1.3}[1.3]{\includegraphics{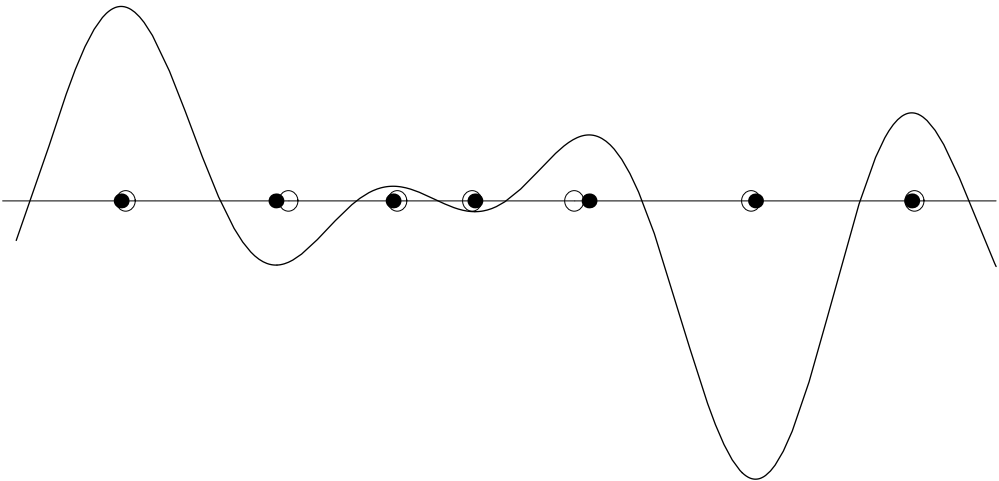}}
\caption{\sf The circles are the midpoints of neighboring zeros of the polynomial~$f$, 
and the dots are the zeros of~$f'$. Think of the circles as a target, which
the zeros of $f'$ often miss.} \label{fig:zfunction}
\end{center}
\end{figure}

It is tempting to think of the zeros of $f'$ as lying close to the midpoint
of neighboring zeros of~$f$.  
That is a useful model, and we analyze the sequence of midpoints and
use it as a tool to study the zeros of the derivative.
A more accurate picture is think of the zeros of~$f'$ as trying to move as far
as possible from the zeros of $f$, and so move toward regions where
there are fewer zeros of~$f$.  See Figure~\ref{fig:zfunction} for an example.  
If $z_j<z_{j+1}$ are consecutive
zeros of $f$, then Rolle's theorem asserts that there will be a zero
$z'$ of $f'$ between $z_j$ and $z_{j+1}$.
If there are a large number of other zeros near $z_{j+1}$,
then $z'$ will be closer to $z_j$.

There are several intuitive reasons for this phenomenon.
Suppose $z_j$ are the zeros of $f$,
listed in increasing order, and consider
\begin{equation} \label{fpf}
\frac{f'}{f}(z)=\sum_k \frac{1}{z-z_k} .
\end{equation}
If $z_j<p<z_{j+1}$ then $f'(p)=0$ if and only if $\frac{f'}{f}(p)=0$.
Substituting $z=p$ into (\ref{fpf}) we see two kinds of terms:
positive for $k\le j$, and negative for $k\ge j+1$.
If there are many zeros near $z_{j+1}$ then all of those negative
terms have to be balanced by something. If there aren't many zeros
near $z_j$, then the only option is
moving $p$ closer to $z_j$.  One can think of  $M=\frac{f'}{f}(\frac12 (z_j+z_{j+1}))$
as a measure of how far the zero of $f'$ will be from the midpoint.
Specifically, 
\begin{equation}
p\approx \begin{cases}
\frac12 (z_j+z_{j+1}) + \frac{M}{8} (z_{j+1}-z_j)^2 & \text{ if } M  \text{ is small} \\
(z_j \, \hbox{ or } \, z_{j+1}) -\frac {1}{M} & \text{ if } M  \text{ is large.}
\end{cases}
\end{equation}

An illustrative example is the polynomial $f$ with zeros at
$\{0,1,\ldots,N\}$.  The first zero of the derivative is very
close to $1/\log N$, and $(f'/f)(\frac12)\approx - \log N$.

A less explicit version of the above explanation uses Gauss' 
electrostatic model.  By (\ref{fpf}), the
zeros of $f'$ can be thought of as points of zero electric 
field for a collection of
equal charges located at~$\{z_k\}$, where like charges repel according
to an inverse linear law.  If there is a region of high concentration
of charges, then these will push the equilibrium points toward
regions of lower charge concentration.

Another approach concerns the relationship between 
the spacing of zeros and the size of the relative maxima and minima
of~$f$.  Roughly speaking, $f$ gets big when there are larger than average
sized gaps between zeros.  But since $f$ is a polynomial with all zeros on a line,
it can't change slope too fast.  This forces the locations of the
relative maxima to shift around so that the slope of $f$  can change appropriately. 
See Figure \ref{fig:zfunction}, particularly near the first local minimum
and the third local maximum.

This approach is interesting because the current view
is that the spacing of zeros is key to understanding the large values
of functions like the Riemann zeta-function.  

In addition to entire functions with zeros on a line, we also 
consider polynomials with all their zeros on a circle.  
We find a close analogy between the two cases.
The motivation here is that much recent work has concerned the 
connection between $L$-functions 
(such as the Riemann $\zeta$-function,
which conjecturally have all their zeros on a line)
and the characteristic polynomials of random
unitary matrices (which have all their zeros on the unit circle).  
The connection between random matrices and $L$-functions is 
discussed extensively in~\cite{cfkrs}.

In the next section we start with an old result that motivates 
the remainder of the paper.  We first consider polynomials, and then
generalize to entire functions of order~1.
Our main results are stated in Section~\ref{sec:longterm}.
In Section~\ref{midpts} we consider a somewhat
simpler averaging process which we later will compare to
differentiation.
In Section~\ref{almost} we  consider the case of zeros which are almost
equally spaced, and determine the rate at which differentiation
evens out the spacing.  
Note that in Sections~\ref{midpts} and~\ref{almost} we treat in 
parallel the cases of zeros on a line and zeros on a circle.  
In Section~\ref{cosineproof} we complete the proof of our main results, and also 
discuss the case that zeros do not lie on a line.

We thank Robin Chapman, Brian Conrey, Nathan Feldman, Chris Hughes, 
and Paul McGuire for helpful conversations,
and the referee for several useful suggestions.

\section{Small gaps become larger}

If $f$ is a function with all its zeros on
a line, define the ``smallest gap between zeros'' function
$$
G[f] = \inf_{j\not=k} |a_j-a_k|
$$
where $a_j$ and $a_k$ are distinct  zeros of $f$, with zeros repeated
according to their multiplicity (so $G[f]=0 $ if $f$ has
a multiple zero).  We will show that
if $f$ is a polynomial with all zeros on a line, then $G[f']\ge G[f]$. 
That is,
the smallest gap between zeros of $f'$ is at least as large as the
smallest gap between zeros of~$f$.  
This is in keeping with the underlying principle that zeros 
of $f'$ tend to ``move away'' from
regions with higher density of zeros of $f$, and to ``move toward''
regions with lower density of zeros of~$f$.  

\subsection{Motivating result}
The following theorem is attributed to M.~Riesz, see~\cite{st}.  
The result was rediscovered several
times~\cite{sz,w1,w2}, including some results which are
quantitatively stronger than the statement below.  
Our proof is similar in spirit, but somewhat
simpler, than previous versions.

\begin{theorem}\label{zpgreater}
Let $f(x)=e^{\alpha x} h(x)$, where $h$
is a  polynomial with only real zeros.
If $\alpha$ and $a$ are real then $G[f'+af]\ge G[f]$,
with strict inequality provided the zeros of $f$ are simple.
\end{theorem}

In particular, zeros of $f^\prime$ are `further apart' than zeros of $f$.

It is implicit, and we give Lagrange's proof following 
equation (\ref{fpoverf}) below, that all
zeros of $f' + a f$ are real.

\begin{proof} We may assume that $f$ has only simple zeros
$a_0<a_1<\cdots<a_J$.
We can write
\begin{equation}
f(x)= e^{\alpha x} \prod_{j=0}^J (x-a_j) 
\end{equation}
and so 
\begin{equation}
\frac{f^\prime + a f}{f}(x) = \alpha + a + 
\sum_{j=0}^J \frac{1}{x-a_j}.
\label{fpoverf}
\end{equation}
We see that all the zeros of $ f^\prime + a f$
are real because if  $x$ has 
positive (negative) imaginary part, then each term in the
sum has negative (positive) imaginary part.
Also, between each 
pair of zeros of $f$ there is exactly one zero of $ f^\prime + a f$, and there
is one additional zero if $a + \alpha \not = 0$.

Let $q<p$ be successive zeros of $f^\prime + a f$.
Assume, for a contradiction,
$p-q\le a_{j+1} - a_j$ for all $j$.
Note that this is equivalent to
\begin{equation}\label{maininequality}
\frac1{p-a_{j+1}}\ge \frac1{q-a_j} ,
\end{equation}
because $p-a_{j+1}$ and $q-a_j$ are either both positive
or both negative, because there is only one zero of $f$ between
$p$ and $q$.

We have
\begin{eqnarray*}
0 &=& \alpha + a+
		\sum_{j=0}^J \frac{1}{p-a_j},\\
\end{eqnarray*}
and
\begin{eqnarray*}
0 &=&  \alpha + a+\sum_{j=0}^J \frac{1}{q-a_j}.
\end{eqnarray*}
Subtracting the two equations gives
\begin{eqnarray}\label{eqn:firstproof}
0 &=& \sum_{j=0}^J \(\frac{1}{p-a_j} - \frac{1}{ q-a_j} \) \cr
&=& \frac1{p-a_0} - \frac1{q-a_J} +
\sum_{j=0}^{J-1} \(\frac{1}{p-a_{j+1}} - \frac{1}{ q-a_j} \) .
\end{eqnarray}
The first two terms above are strictly positive, 
and by (\ref{maininequality}) all the
terms in the sum are nonnegative.  So the right side is strictly 
positive, which is a contradiction.
\end{proof}

Repeatedly applying the result leads to various linear combinations
of $f$ and its derivatives which satisfy the corresponding inequality.
For example, $G[f-cf'']\ge G[f]$
provided $c\ge 0$.

\subsection{Entire functions of order 1}

We wish to generalize the above result to functions having infinitely
many zeros.  The appropriate class to consider is entire functions
of order at most~1.  We recall the standard terminology and properties
of these functions.
For more details, see~\cite{boas}.

If $f$ is an entire function, then $f$~\emph{is of order $\rho$} if
\begin{equation}
\limsup_{r\to\infty} \frac{\log\log M(r)}{\log r} =\rho ,
\end{equation}
where $M(r)=\max_{|z|=r} |f(z)|$. That is, $f(z)=O(\exp(|z|^{\rho+\delta}))$
for all $\delta>0$ and no $\delta<0$.

 If $f$ is of order $\rho$, then
$f$ \emph{has type $\tau$} if
\begin{equation}
\limsup_{r\to\infty} r^{-\rho} \log M(r) =\tau .
\end{equation}
That is, $f(z)=O(\exp((\tau+\delta) |z|^{\rho}))$
for all $\delta>0$ and no $\delta<0$. 
We say that $f$ has minimal (respectively, finite or maximal) type
if $\tau=0$ (respectively, $0<\tau<\infty$ or $\tau=\infty$).

For example, $\cos(z)=\frac12 (e^{iz}+e^{-iz})$, so $\cos(z)$ has
order~1 and type~1.  And by Stirling's formula, 
$1/\Gamma(z)$ has order~1 and maximal type.

Of particular interest is the connection between order, type, and
the distribution of zeros.  Let $(z_n)$ denote the zeros
of~$f$, with $z_j\le z_{j+1}$ and zeros repeated 
according to their multiplicity.
Denote by $n(r)$ the number of zeros of $f$ in $|z|\le r$.  If $f$
has only real zeros, then let $n_+(r)$ (respectively, $n_-(r)$)
denote the number of zeros in $(0,r]$ (respectively, $[-r,0)$).
A characterization of functions of finite type is due to
Lindel\"of:
\begin{proposition}(Lindel\"of)  If $\rho$ is an positive integer, then the
entire function $f$ of order $\rho$ is of finite type if and
only if $n(r)=O(r^\rho)$ and the sum over zeros
\begin{equation}
S(r)=\sum_{|z_n|\le r} z_n^{-\rho}
\end{equation}
is bounded.
\end{proposition}

That $\cos(z)$ has finite type and $1/\Gamma(z)$ has maximal
type follows directly from Lindel\"of's theorem.

The Riemann $\zeta$-function  is not entire because it 
has a simple pole at $s=1$, but $(s-1)\zeta(s)$ is entire.  
From the Dirichlet
series representation $\zeta(s)=\sum_{n=1}^\infty n^{-s}$ it follows 
that the $\zeta$-function
is bounded for $\Re(s)>1+\delta$.  By the functional equation
$$
\zeta(s)=2^s\pi^{s-1}\sin(\tfrac12 s\pi) \Gamma(1-s)\zeta(1-s)
$$
and the fact that the sine and $\Gamma$-function have order~1, it follows that
the $\zeta$-function  has order~1.  

The Riemann $\zeta$-function has maximal type because $\zeta(s)$ has 
zeros at the negative even integers and 
$\frac{1}{2\pi} T\log(T/2\pi e)+O(\log T)$ zeros with imaginary
part $0<t<T$ and real part $0<\sigma<1$, 
so it fails both tests in Lindel\"of's theorem.

The Riemann $\xi$-function is given by
\begin{equation}
\xi(s)= \tfrac12 s(s-1)\pi^{-\frac12 s} \Gamma(\tfrac12 s) \zeta(s).
\end{equation}
It satisfies the functional equation $\xi(s)=\xi(1-s)$.
The $\xi$-function is entire of order~1 and has the same non-real zeros
as the $\zeta$-function, so it has maximal type.  The Riemann Hypothesis
is the assertion that all zeros of the $\xi$-function have real
part~$\frac12$.  

Sometimes it is convenient to use the Riemann
$\Xi$-function, given by $\Xi(z)=\xi(\frac12+iz)$.  The $\Xi$-function
is an even function which is real on the real axis, 
is entire of order~1 and maximal type.  Assuming the Riemann Hypothesis,
$\Xi$ has only real zeros and 
$n_+(r)=n_-(r)=({r}/{2\pi})\log(r/2\pi e)+O(\log r)$.

In the generalization of Theorem~\ref{zpgreater} to entire functions of order 1
we will make use of the Hadamard product.  Suppose $f$ is entire of order~$\rho$
and let $N=[\rho]$ be the largest integer less then or equal to~$\rho$.
Then there exists a polynomial~$Q$ of degree at most~$N$ such that
\begin{equation}
f(z)=e^{Q(z)}z^{n_0} \prod_n \(1-\frac{z}{z_n}\) \exp\(P_N\(\frac{z}{z_n}\)\),
\end{equation}
where $P_N(z)=z^{}+\frac12 z^2+\cdots+\frac1N  z^{N}$ and $n_0$ is the order of the zero
of $f(z)$ at $z=0$.
 
\subsection{Generalization of Theorem~\ref{zpgreater} to entire functions of order 1}\label{generalization}
We have

\begin{theorem}\label{thm:order1version}
Suppose $f$ is an entire function of order~1 which is real on the
real axis and has only real zeros.  If $q<p$ are consecutive
zeros of~$f' + a f $, $a\in \R$, then
\begin{equation}
\inf z_{n+1}-z_n \le p-q \le \sup z_{n+1}-z_n .
\end{equation}
If the zeros of $f$ are simple and equality holds for one of the 
inequalities, then it also holds for
the other, and $f(z)=A e^{B z} \cos(C z + D)$ for some real 
$A$, $B$, $C$, and~$D$.
\end{theorem}

So, if the zeros of $f$ are not equally spaced, then differentiation makes
the smallest gaps larger and the largest gaps smaller.  
It might not be obvious that if the zeros of $f$ are
equally spaced then the same holds for~$f'$, but one can check that
the set of functions of the form $A e^{B z} \cos(C z + D)$ is closed
under differentiation.  In particular, the factor $e^{Bx}$ shifts all of the 
zeros of the derivative of $\cos(C x + D)$ by the same amount.

\begin{proof} 
We may assume that $f(0)\not=0$ and the zeros of $f$ are simple.
We have the Hadamard factorization
\begin{equation}\label{eqn:order1hp}
f(z)=Ae^{Bz} \prod_n \(1-\frac{z}{z_n}\) \exp\(\frac{z}{z_n}\),
\end{equation}
so
\begin{equation}
\frac{f'+af}{f}(z)=a+B+\sum_n \frac{1}{z-z_n}-\frac{1}{z_n} .
\end{equation}
Imitating the previous proof, the terms $a$, $B$, and $z_n$ cancel
and equation~(\ref{eqn:firstproof}) becomes
\begin{equation}\label{zerosum}
0=\sum_j \frac{1}{p-z_{j+1}}-\frac{1}{q-z_j} .
\end{equation} 
By~(\ref{maininequality}) every term in the sum is nonnegative,
and we have a contradiction unless $p-q=z_{j+1}-z_j$ for all~$j$.
This proves the left hand inequality in the Theorem.  The proof of
the other inequality is identical.

Finally, if $p-q=z_{j+1}-z_j$ for all~$j$ then the zeros of $f$ are
equally spaced, so we recognize the Hadamard product of $f$ 
to be of the stated form.
\end{proof}

A shortcoming of Theorem~\ref{thm:order1version} is that the inequalities 
can be vacuous.  For example, for 
the Riemann~$\Xi$-function $\inf z_{n+1}-z_n=0$.
To formulate a ``local'' version of Theorem~\ref{thm:order1version}
we must be able to speak of the ``local density'' of zeros.
Suppose
\begin{equation}\label{ndensity}
n_{\pm}(r) =  L(r) + E(r), 
\end{equation}
where $L(r)$ is a nice function and $E(r)=o( L(r) )$ should be thought
of as an error term.  Then
$L'(r)$ is the density of zeros near~$r$, meaning that 
$L'(r)^{-1}$ is the average gap between zeros near~$r$.
We want to show that if $f'$ has zeros which are closer together
than the average spacing, then in some neighborhood 
$f$ will also have zeros which 
are closer together than average.  The following Theorem 
asserts that the neighborhood  can be taken to be of size  $E(r)/L'(r)$.

\begin{theorem}\label{thm:localversion}
Suppose $f$ is an entire function of order~1 which is real on the
real axis, has only real zeros, and $n_{\pm}(r)$ satisfies
(\ref{ndensity}) 
where $L$ is an increasing smooth function
such that $ L^{(j)} (r)\ll r^{-j}L(r)$.
Suppose $\kappa<K<1$ and 
$q<p$ are consecutive
zeros of~$f'+a f$, $a\in \R$,  such that $p-q<\kappa L'(q)^{-1}$. Then
there exist zeros~$z_n$,~$z_{n+1}$ of $f$ such that
$|z_n-q|\ll E(q)/L'(q)$ and 
$z_{n+1}-z_n < K  L'(q)^{-1}$.
\end{theorem}

For the
Riemann~$\Xi$-function we have 
$L(r) =(2\pi)^{-1} r\log (r/2\pi e)$ and $E(r)=\log r$,
giving this
\begin{corollary}\label{logj}
Suppose $0<\kappa<K<1$ and $\gamma_1'<\gamma_2'$ are 
consecutive zeros of $\Xi^{(n)}$
with $\gamma_2'-\gamma_1'<2\pi \kappa/\log\gamma_1$.  Then there exist consecutive
zeros $\gamma_1<\gamma_2$ of $\Xi$ with $|\gamma_1-\gamma_1'|\ll 1$ such that
$\gamma_2-\gamma_1<2\pi K/\log\gamma $.
\end{corollary}

The Corollary is motivated by the problem of Landau-Siegel zeros.
Conrey and Iwaniec~\cite{ci} have shown that
if there exists $\kappa <\frac12$ and $\delta>0$ 
such that $\gamma_2-\gamma_1 <2\pi \kappa /\log\gamma$
for $\gg T \log^\delta T$ pairs of zeros $0<\gamma_2,\gamma_1<T$ of the 
$\Xi$-function,
then one can find an effective lower bound for the class number of the
imaginary quadratic field~$\Q(\sqrt{-d})$.  By the Corollary, 
if we had $\gg T\log^\delta T$ pairs of zeros 
$\gamma_2'-\gamma_1' <2\pi K/\log\gamma$
for $\Xi^{(n)}$, for any fixed $n>0$, then we would have $\gg T$
such pairs of zeros 
for~$\Xi$.  This falls short of the requirements for 
Conrey and Iwaniec's result.  However, it may be possible 
to modify the proof of Theorem~\ref{thm:localversion} to show that
either $\Xi$ has a pair of zeros which are closely separated and
very near to~$q$, or else $\Xi$ has a large number of 
zeros which are not too far from~$q$.

\begin{proof}[Proof of Theorem~\ref{thm:localversion}]  
For a contradiction, suppose $z_{j+1}-z_j > p-q  $ for $|q-z_j|<X$.
We can reuse everything in the previous proof up to~(\ref{zerosum}):
\begin{eqnarray}
0 &=& \sum_{j} \(\frac{1}{p-z_{j+1}} - \frac{1}{ q-z_j} \) \cr
&=& \(\sum_{|q-z_j|\le X} + \sum_{q-z_j<- X} + \sum_{z_j-q> X}\)  
			\(\frac{1}{p-z_{j+1}} - \frac{1}{ q-z_j} \) \cr
&=&S_1+S_2+S_3,
\end{eqnarray}
say. 

By assumption, $S_1\ge 0 $.
We now estimate $S_3$, the treatment of $S_2$ being identical. Note that
\begin{eqnarray}\label{telescope}
\frac {1}{p-z_{j+1}} -\frac{1}{q-z_j} &=&
\frac {1}{q-z_{j+1}} -\frac{1}{q-z_j} +\frac {p-q}{(p-z_{j+1})({q-z_{j+1}})} \cr
&=& \frac {1}{q-z_{j+1}} -\frac{1}{q-z_j} +\frac{p-q}{(q-z_{j+1})^2}
+O\(\frac{(p-q)^2}{(q-z_{j+1})^2}\) .
\end{eqnarray}
Let $Z_J$ denote the smallest zero of $f$ in~$[q+X,\infty)$, so the
sum in $S_3$ begins with~$z_J$.
In the sum $S_3$, the first two terms in~(\ref{telescope}) telescope, 
giving
\begin{equation}\label{presteiltjes}
S_2 = \frac{1}{z_J-q} + (p-q)(1+O(L'(q)^{-1})\sum_{z_j-q> X} \frac{1}{(q-z_{j+1})^2} .
\end{equation}
Writing the sum as a Steiltjes integral and integrating by parts we have
\begin{eqnarray}
\sum_{z_j-q> X} \frac{1}{(q-z_{j+1})^2} &=&
\int_{z_{J+1}}^\infty \frac{dn_+(r)}{(r-q)^2} \cr
&=& \int_{z_{J+1}}^\infty \frac{L'(r)}{(r-q)^2}\, dr 
		+ \int_{z_{J+1}}^\infty \frac{dE(r)}{(r-q)^2} \cr
&=& \frac{L'(z_{J+1})}{z_{J+1}-q} 
                + \frac{E(q)}{(z_{J+1}-q)^2} 
+ \int_{z_{J+1}}^\infty \frac{L''(r)}{r-q}\, dr 
		+ 2 \int_{z_{J+1}}^\infty \frac{E(r)}{(r-q)^3} \,dr \cr
&=& (1+o(1))\frac{L'(z_{J+1})}{z_{J+1}-q} 
                + \frac{E(z_{J+1})}{(z_{J+1}-q)^2} .
\end{eqnarray}
In the last step we used $E(r)=o(L(r))$ and $L''(r)\ll L'(r)/r$.

Adding all the terms  and using the fact that $L'(z_{J+1})\sim L'(q)$ and
$E(z_{J+1})\sim E(q)$ we have
\begin{equation}
S_3 \sim  \frac{1}{z_{J}-q} - \frac{(q-p)L'(q)}{z_{J+1}-q} 
         - \frac{E(q) (q-p) (z_{J+1}-q)^{-1}}{z_{J+1}-q} .
\end{equation}
We assumed $q-p\le \kappa L'(q)^{-1}$, so the first two terms 
sum to something larger than $\delta/(z_{J+1}-q)$ for some $\delta>0$.  
The numerator of the last term is
bounded by 
$\kappa E(q) L'(q)^{-1} X^{-1} $ which is smaller than $\delta$ provided
$ X>  \delta^{-1} \kappa E(q) L'(q)^{-1}$.  For such $X$ we have $S_3>0$,
which gives a contradiction.  We have $X\ll  E(q)/L'(q)$, as claimed.
\end{proof}

\subsection{Long-term behavior}\label{sec:longterm}
As described in the previous section, differentiation causes the
small gaps between zeros of $f$ to become larger and the large gaps to
become smaller.  Thus, one might expect that under repeated differentiation
the zeros of $f^{(n)}$ will approach equal spacing, which would equal
the average spacing between zeros of~$f$.  And if the zeros of $f^{(n)}$
are approaching equal spacing then $f^{(n)}(z)$ should be approaching
$Ae^{Bz}\cos(Cz+D)$, where $A$, $B$, $C$, and $D$ may depend on~$n$. 

In order for the above argument to hold, 
there must be an appropriate sense of average spacing between
zeros of $f$.  Below we show that it is sufficient to
have $n_+(r)\sim n_-(r)\sim \kappa r$
for some~$\kappa>0$, but it is possible that this condition
can be weakened.  Since $Ae^{Bz}\cos(Cz+D)$ has finite type,
one might expect that this would also be a necessary condition
on~$f$, but this is not the case.  The assumptions in the following result
are stronger than the first condition in Lindel\"of's theorem, but
weaker than the second.

\begin{theorem}\label{cosinetheorem}
Suppose $f$ is an entire function of order~1 which is real on the
real axis, has only real zeros, and
$n_+(r)\sim n_-(r)\sim \kappa r$.
Then there exist sequences $(A_n)$, $(B_n)$, and $(D_n)$ with
$D_n$ bounded, such that
\begin{equation}
\lim_{n\to\infty} A_n e^{B_n z} f^{(n)}(\kappa^{-1} z+D_n) = \cos(\pi z) ,
\end{equation}
uniformly for $|z|\le X$ for any fixed $X>0$.  In particular, the
zeros of $f^{(n)}$ approach equal spacing.
\end{theorem}

We give the proof in Section~\ref{cosineproof}.

It seems that $A_n$ in the theorem is increasing or decreasing
according to whether $n(r)-2\kappa r$ is generally negative
or positive, but we have not succeeded in finding a precise
statement.  
The Theorem also suggests that the regularity in the distribution
of the zeros of $f$ should lead to a great regularity in the Taylor
series coefficients of~$f$.  This seems worth exploring further.

A trivial example of the theorem is
when
$f(z)=A e^{Bz} \cos(z+D)$ for some $A,B,D$.
We have $f'(z)=\sqrt{B^2+1} f(z+\varphi)$ where
$\cos(\varphi)=B/\sqrt{B^2+1}$ and $\sin(\varphi)=1/\sqrt{B^2+1}$.
So the conclusion of the theorem holds with
$A_n=A^{-1}(B^2+1)^{-\frac{n}{2}}$ and $D_n \equiv -D-n \varphi \ (\bmod 2\pi)$.
Somewhat more generally, suppose
$f(z)=A e^{Bz}\( \cos(z+D)+E\)$ where $|E|\le1$.  This function has zeros
consisting of two interlaced sequences which separately are equally spaced.
The conclusion of the theorem holds with the same $A_n$ and~$D_n$ as in
the previous example.
Similarly, if $f(z)=F(z)\cos(z)$ where $F$ is a degree $N$ polynomial,
then the conclusion holds with $A_n=n^{-N}$ and 
$D_n\in\{-\frac{\pi}{2},0,\frac{\pi}{2},\pi\}$.

A less trivial example is the Bessel function~$J_n(z)$.  We
have $n_{\pm}(r)=\pi^{-1} r +O(1)$, so the theorem applies.  
In this case the theorem can be verified
directly because
$$
J_n(\pi z)=\frac{1}{\pi} \int_0^\pi
\cos(\pi z\sin t - n t)\, dt ,
$$
so for an even number of derivatives
\begin{eqnarray}\label{besselcosine}
J_n^{(2k)}(\pi z)&=&\frac{(-\pi)^k}{\pi} \int_0^\pi
\sin^{2k}(t)
\cos(\pi z\sin t - n t)\, dt \cr
&\sim& (-\pi)^k \sqrt{\frac{\pi}{k}} \cos\(\pi z - \frac{n\pi}{2} \),
\end{eqnarray}
because $\sqrt{k/2\pi}\sin^k t\, dt$ is approaching a
Dirac $\delta$-function at $t=\pi/2$.  So the conclusion of
the theorem holds with $A_{2k} = \pi^{-k}\sqrt{k/\pi}$ and
$D_{2k}\equiv \frac{n\pi}{2}\mod 2\pi$, and similarly for $k$~odd.
Equation~(\ref{besselcosine}) can also be obtained from the
Taylor series for~$J_n$.

Theorem~\ref{cosinetheorem} asserts that the zeros of~$f^{(n)}$
approach equal spacing.  The approach to equal spacing is actually quite fast.
The general case is somewhat intractable, so we will assume
that the zeros are close to equally spaced and do a first
order approximation.  We also
suppose further that $f$ is approximately an odd function, although this is
not essential and we describe the necessary modifications during the proof.
Thus,
\begin{equation}
\label{feps}
f(z) = C\,(z-\ve_0) \prod_{j\not=0} \(1- \frac{z}{j +\ve_j}\) \exp\(\frac{z}{j+\ve_j}\),
\end{equation}
where $|\ve_j| < \ve$ for some small $\ve$.

\begin{theorem}\label{oneoverj}
Suppose  $f$ is given in (\ref{feps}) and let
$z_n^j$ be the zeros of the $j$th derivative $f^{(j)}$.
Then $z^j_{n+1}-z^j_n=1+O(\ve j^{-1})$.
\end{theorem}

We give the proof in Section~\ref{rate}.

If the zeros of $f$ do not have a nonzero average spacing,
then the results in this section do not apply.
Suppose, for example,  $n_+(r)$ and $n_-(r)$ grow faster than
linearly.
The ``overall'' average spacing  is~$0$,
and under repeated differentiation there is a competition 
between zeros moving toward the origin
and zeros trying to become locally equally spaced.
If $n_+(r)$ and $n_-(r)$ are sufficiently nice then
there is the possibility that an analogue of 
Theorem~\ref{cosinetheorem} may hold.  This was established
by Haseo Ki~\cite{ki2} in the special case of the
Riemann $\Xi$-function:

\begin{theorem}(Haseo Ki~\cite{ki2}) There exist
sequences $A_n$ and  $C_n$, with $C_n \to 0$ slowly,
such that
\begin{equation}
\lim_{n\to\infty} A_n\, \Xi^{(2n)}(C_n z)
=
\cos(z) 
\end{equation}
uniformly on compact subsets of $\mathbb C$.
\end{theorem}

The above result was conjectured in an earlier verision of
this paper.  The $\Xi$-function has
$n_\pm(r)=(2\pi)^{-1} r\log(r/2\pi e)+O(\log r)$ so the
zero counting function is particularly well-behaved.  However, Ki's proof
makes use of the fact that the $\Xi$-function
has a nice representation as a Fourier transform, rather
than directly using properties of the zeros.
The theorem suggests that $\log(\xi^{(n)}(\frac 12))$ should grow very
regularly and not too much faster than linearly
as $n\to\infty$.  Rick Kreminski~\cite{kr1,kr2} has calculated the first
$490$ derivatives and his data appears to grow approximately like~$n\log n$.

More generally we conjecture:

\begin{conjecture}\label{Xiconjecture} Suppose $f$ is a real entire
function of order~1 having only real zeros, with  $n_+(r)$ and $n_-(r)$
sufficiently nice and $n_+(r)-n_-(r)$ not too large. Then there
exist sequences $A_n$, $B_n$, $C_n$, and $D_n$ with $D_n/C_n$ bounded,
such that
\begin{equation}
\lim_{n\to\infty} A_n\, e^{B_n z} f^{(2n)}(C_n z + D_n)
=
\cos(z) 
\end{equation}
uniformly on compact subsets of $\mathbb C$.
\end{conjecture}

\section{Averaging, instead of differentiating}\label{midpts}

Differentiation is a process which takes one sequence of points
(the zeros of $f$) and replaces it with another sequence of points
(zeros of $f'$), such that the two sequences interlace each other.
For comparison, we also consider the much simpler process
of making a new sequence from the midpoints of neighboring elements
of the given sequence.

\subsection{Averaging on the line}

Suppose $(x_n)$
is an increasing sequence.  Form new sequences $(x^{j}_n)$ where
$x^{0}_n=x_n$ and
$x^{j+1}_n=\frac12 \(x^{j}_n+x^{j}_{n+(-1)^j}\)$. 
That is, the terms of each new sequence are the midpoints
of consecutive terms in the previous sequence.
The subscript ${n+(-1)^j}$ is designed so that if 
 $(x^{j}_n)$ is equally spaced, then $x^{j+2}_n=x^{j}_n$.
If that subscript was ${n+1}$ then the equal spaced case
would give $x^{j+2}_n=x^{j}_{n-1}$, and the terms would
be ``drifting to the right.''

\begin{theorem}\label{sqrtj}
Suppose $(x_n)_{n\in \Z}$
is a sequence with $x_n=n+\ve_n$, where
$\ve_n \ll E(|n|)$, as $n\to\pm\infty$, for some increasing function~$E$.
If $E(n)=o(n)$, that is $x_n\sim n$, then
as $j\to\infty$ the $j^{th}$ midpoint sequence $(x^{j}_n)$
approaches equal spacing.
Furthermore, if $E(n) \ll n^\theta$ with $0\le \theta<1$
then
$$
x^{j}_{n+1} - x^{j}_n = 
1+O\((n^\theta+j^{\frac12 \theta})j^{-\frac12}\) .
$$ 
\end{theorem}

\begin{proof}
It is straightforward to show by induction that
$$
x^{j}_n = 2^{-j} \sum_{0\le k\le j}
\binom{j}{k} x_{n+[j/2]-k}.
$$
Thus,
\begin{eqnarray}\label{binomsum}
x_{n+1}^j-x_n^j &=& 2^{-j} \sum_{0\le k\le j}
\binom{j}{k} \(x_{n+1+[j/2]-k} - x_{n+[j/2]-k}\) \cr
&=& 2^{-j} \sum_{0\le k\le j}
\binom{j}{k} 
+2^{-j} \sum_{0\le k\le j} \binom{j}{k} \(\ve_{n+1+[j/2]-k} - \ve_{n+[j/2]-k}\) \cr
&=& 1
+2^{-j} \sum_{0\le k\le j} \ve_{n+1+[j/2]-k}\left(\binom{j}{k}- \binom{j}{k-1} \) \cr
&&
\mathstrut + 2^{-j} \(
\ve_{n+[j/2]-j} - \ve_{n+[j/2]}\) \biggr. \cr
&=& 1
+2^{-j}\sum_{-\frac{j}2 \le \ell \le \frac{j}2 } \ve_{n+\ell }\,
\binom{j}{\frac{j}{2}+\ell }
\frac{-\ell+1}{\frac{j}{2}-\ell+1}\ 
+O(2^{-j}(n+j)).
\end{eqnarray}
We may assume that $E$ is an increasing function and
$E(a+b)\ll E(a)+E(b)$ for $a,b>0$.  So the
sum in~(\ref{binomsum}) is bounded by
\begin{equation}
\sum_{0 \le \ell \le \frac{j}2 } \(E(n)+E(\ell)\)\,
\binom{j}{\frac{j}{2}+\ell }
\frac{\ell}{\frac{j}{2}-\ell+1} .
\end{equation}

Now use fact that the binomial distribution approaches the
Gaussian to estimate the binomial coefficient:
\begin{equation}
\binom{j}{k} \sim \frac{2^{j+\frac12}}{\sqrt{\mathstrut \pi j}} \,
		\exp\(-\frac{2(k-\frac{j}{2})^2}{j} \).
\end{equation}
So,
\begin{eqnarray}\label{eqn:finalbonomialbound}
\sum_{0 \le \ell \le \frac{j}2 } E(\ell)
\binom{j}{\frac{j}{2}+\ell }
\frac{\ell}{\frac{j}{2}-\ell+1}
&\ll&
 \frac{2^{j}}{j^{\frac12}} 
\sum_{0 \le \ell \le \frac{j}2 } E(\ell)
                \exp\(-\frac{2\ell^2}{j}\)
\frac{\ell}{\frac{j}{2}-\ell+1}\cr
&\ll&  \frac{2^{j}}{j^{\frac12}} 
\int_0^\frac{j}{2} 
E(t)
                \exp\(-\frac{2t^2}{j}\)
\frac{t}{\frac{j}{2}-t+1}\, dt \cr
&=& 2^{j} j^{\frac12}
\int_0^\frac{\sqrt{\mathstrut j}}{2} 
E(t\sqrt{\mathstrut j})
                \exp\(-{2t^2}\)
\frac{t}{\frac{j}{2}-t\sqrt{\mathstrut j}+1}\, dt . 
\end{eqnarray}
We may assume $E(ab)\ll E(a)E(b)$ for $a,b>1$, 
so~(\ref{eqn:finalbonomialbound}) 
is bounded by
\begin{eqnarray}
2^{j} j^{\frac12} E(\sqrt{\mathstrut j})
\int_0^\frac{\sqrt{\mathstrut j}}{2} 
                \exp\(-{2t^2}\)
\frac{t^2}{\frac{j}{2}-t\sqrt{\mathstrut j}+1}\, dt 
\ll 2^{j} j^{-\frac12} E(\sqrt{\mathstrut j}) ,
\end{eqnarray}
the integral being seen to be $\ll j^{-1}$ by breaking it at $j^{\delta}$,
for any $0<\delta<\frac14$, and 
estimate the two parts separately.

To finish the proof, note that if $E(n)=o(n)$ then 
$j^{-\frac12} E(\sqrt{\mathstrut j})=o(1)$, and if
$E(n)\ll n^\theta$ then the other estimates follows immediately.
\end{proof}

If $x_n=n+O(1)$ then the above shows that the $j$th midpoint sequence 
approaches equal spacing with a discrepancy of order $j^{-\frac12}$.
This convergence to equal spacing is actually quite slow, and in some sense
it is the slowest way to ``even out'' a sequence.  
In particular,  Theorem~\ref{oneoverj} says
taking successive derivatives evens out the sequence much faster,
with a discrepancy of order  $j^{-1}$.  These topics are discussed in 
Section~\ref{rate}.

\subsection{Averaging on the circle}\label{circle}

As mentioned in the Introduction, 
it is believed that there is a close analogy 
between $L$-functions, which (conjecturally) have all their zeros on a line, and the
characteristic polynomials of random unitary matrices, which have all
their zeros on a circle.  See~\cite{cfkrs}.
It is obvious that the analogy must break 
down if pushed too far, for a polynomial has only finitely many zeros.
In this section we consider the averaging of points on a circle,
in analogy to the averaging on a line in the previous section.
Here we find that the discrepancy from equal spacing vanishes
exponentially, in sharp contrast to the~$j^{-\frac12}$ of the linear case.
This exponential convergence is a general property of iterative averaging
procedures on finite sets, as described in Section~4.4 of~\cite{davis}.

Suppose $p_1,\ldots,p_\ell$ are points on the unit circle, and form new
sets of points on the circle  $(p^j_n)$ in analogy to the previous section, where
$0\le p^j_n<2\pi$ is interpreted as an angle.
If $j<n$ then this process is indistinguishable from the averaging
process in the previous section.  But when $j>n$ the fact that 
there are only $n$ points on a circle comes into play.  This causes the averaging
to even out the spacing more rapidly.

\begin{proposition}  Let  $(p_n)$ and  $(p^j_n)$ be sequences of $\ell$ points
on the unit circle as described above.  Then
$p_{n+1}^{j}-p_{n}^{j} = \frac{2\pi}{\ell} + O((1-\frac{\pi^2}{2\ell^2})^j)$.  
\end{proposition}

\begin{proof}
Starting exactly as before, but using the fact that $p_j=p_{k}$ if 
$j\equiv k\bmod \ell$, we have
\begin{eqnarray*}
p_{n+1}^j-p_n^j &=& 2^{-j} \sum_{0\le k\le j}
\binom{j}{k} \(p_{n+1+[j/2]-k} - p_{n+[j/2]-k}\) \\
&=& 2^{-j} \sum_{m=1}^\ell (p_{m+1}-p_m)
\sum_{\ontop{0\le k\le j}{k\equiv M \bmod \ell}}
\binom{j}{k} ,
\end{eqnarray*}
where $M=n+1+[j/2]-m$.  

Now let $\rho$ be a primitive $\ell^{th}$ root of~$1$, and let
$\rho_j=\rho^j$ for $0\le j\le \ell-1$.
Using the binomial theorem and the
fact that $\sum_k \rho_k^n =0$ except 
when $n\equiv 0\bmod \ell$, we have
\begin{eqnarray}
\ell \sum_{\ontop{0\le k\le j}{k\equiv M \bmod \ell}}
\binom{j}{k} 
&=&
(1+1)^j + \rho_{M_1}(1+\rho_1)^j+\cdots+\rho_{M_{\ell-1}}(1+\rho_{\ell-1})^j \\
&=&
2^j + O\(2^j \ell \(1-\frac {\pi^2}{2\ell^2}\)^j \),
\end{eqnarray}
where $\rho_{M_j}$ is an ordering of the $\rho_k$ depending on~$M$.
The last step follows from $|1+e^{2 \pi i/\ell}| = 2(1-\frac {\pi^2}{2\ell^2}) + O(\ell^{-4})$.

Finally, use the fact that $\sum_{m=1}^\ell (p_{m+1}-p_m)=2\pi$ to 
finish the proof.
\end{proof}

\section{Almost equally spaced zeros}\label{almost}

Suppose $f(z)$ is an entire function of order~1 which is real on the
real axis, has only real zeros, and the zeros have average spacing~1.  
As we repeatedly differentiate, the zeros will approach equal spacing,
and we want to determine the rate at which this occurs.
We will see that the approach to equal spacing is much faster than
in the midpoint process of Section~\ref{midpts}.

We treat in detail the case that the zeros of $f$ are close to
equally spaced, both on the line and on the circle.  At the end of this 
section we determine the rate at which differentiation evens
out zero spacing.

\subsection{Almost equally spaced zeros on a line}\label{almostline}
Suppose
\begin{equation}
f(z) = C\,(z-\ve_0) \prod_{j\not=0} \(1- \frac{z}{j +\ve_j}\) \exp\(\frac{z}{j+\ve_j}\),
\end{equation}
where $|\ve_j| < \ve$ for some small $\ve$.
That is, $f$ is approximately an odd function with zeros close to
equally spaced.

\begin{theorem}\label{alphak}
In the notation above, $f'$ has zeros at $k+\frac{1}{2}+\alpha_k$,
where
$$
 \alpha_k = \frac{4}{\pi^2} \sum_{j} \frac{\ve_{j+k}}{(2j-1)^2} +O(\ve^2), 
$$
and $f''$ has zeros at $k+\beta_k$, where
$$
\beta_k=\frac13 \ve_k +\frac{2}{\pi^2} \sum_{j\not = 0} \frac{\ve_{j+k}}{j^2} +O(\ve^2).
$$
\end{theorem} 

In the proof of Theorem~\ref{alphak} we will require the following standard
formulas:
\begin{equation}\label{standardformulas}
\sum_j \frac{1}{(2j+1)^2} = \frac{\pi^2}{4},
\  \ \ \ \ \  \hbox{and} \ \ \ \ \ \ \ \ 
\sum_j \frac{1}{(2j+1)^2 (2(j+n)+1)^2} = 
\begin{cases}
\frac{\pi^2}{8 n^2} & n\not=0 \cr
\frac{\pi^4}{48} & n=0 .
\end{cases}
\end{equation}

\begin{proof}[Proof of Theorem~\ref{alphak}]
We have
$$
\frac{f'}{f}(z) ={\sum_{j}}' \frac{1}{z - j -
  \ve_j} + \frac{1}{j+\ve_j} .
$$
Here and below, $\sum'$ means that the $j=0$ term must be modified
in the obvious way.
The zeros of $f'(z)$ will be approximately halfway between the 
zeros of $f$, so suppose one of the zeros is $k+\frac12+\alpha_k$.

Note:  if we didn't assume that $f$ was odd, and so a factor 
of $e^{a z}$ occurs in our function, then the zeros of $f'$ would not be 
halfway between the zeros of $f$.  But the zeros of $f'$ would all be
shifted from the midpoint by the same amount, 
and all of the calculations below would 
work with a slight modification.
So,
\begin{eqnarray}
0 = \frac{f'}{f}(k+\frac12 + \alpha_k) &=& {\sum_{j}}'
\frac{1}{k+\frac12 +\alpha_k - j - \ve_j} + \frac{1}{j+\ve_j}\cr
&=& {\sum_{j}}'
\frac{2}{1+2(k-j)} + \frac{4(\ve_j-\alpha_k)}{(1+2(k-j))^2} + \frac{1}{j+\ve_j}
+O(\ve^2+\alpha_k^2)\cr
&=& \sum_{j} \frac{4(\ve_j-\alpha_k)}{(1+2(k-j))^2} +O(\ve^2+\alpha_k^2)\cr
&=& 4\sum_{j} \frac{\ve_{j+k}}{(2j-1)^2} 
- \pi^2 \alpha_k 
+O(\ve^2+\alpha_k^2).
\end{eqnarray}
Solving for $\alpha_k$ gives the first formula in the Theorem.
The $O(\ve^2)$ error term follows because the above formula shows that
$\alpha_k=O(\ve+\ve^2+\alpha_k^2)$, but 
$\ve^2 =o(\ve)$ and $\alpha_k^2=o(\alpha_k)$, so $\alpha_k=O(\ve)$.

To understand the effect of the second derivative on the
zeros, we must iterate the above formula.
The second derivative $f''$ will have zeros near the integers, so
suppose there is a zero at $k+\beta_k$.  By the first formula in the Theorem,
and using formulas~(\ref{standardformulas}),
\begin{eqnarray}
\beta_k &=& \frac{4}{\pi^2} \sum_{j} \frac{\alpha_{j+k-1}}{(2j-1)^2} +O(\ve^2)\cr
&=& \frac{4}{\pi^2} \sum_{j} \frac{1}{(2j-1)^2}
 \frac{4}{\pi^2} \sum_{m} \frac{\ve_{m+j+k-1}}{(2m-1)^2} +O(\ve^2) \cr
&=&  \frac{16}{\pi^4} \sum_{n} \ve_{n+k}  
\sum_{j} \frac{1}{(2j-1)^2}  \frac{1}{(2(n-j+1)-1)^2} +O(\ve^2) \cr
&=&\frac13 \ve_k +\frac{2}{\pi^2} \sum_{n\not = 0} \frac{\ve_{n+k}}{n^2} +O(\ve^2).
\end{eqnarray}
as claimed.

\end{proof}

\subsection{Almost equally spaced zeros on a circle}\label{almostcircle}

If the polynomial $f$ has all its zeros on a circle, then
$f'$ has all its zeros strictly inside the circle, except at points
where $f$ has a multiple zero.  This follows from the famous Gauss-Lucas
theorem that the zeros of $f'$ lie inside the convex hull of the zeros
of~$f$.  Thus, if we are to find an analogy with the case of zeros on
a line, we must do something slightly different than differentiation.

\begin{lemma}
If $f$ is a degree~$n$ polynomial with all its zeros on the unit circle,
then
$$
f^{j}(z):=\(z \frac{d}{dz}\)^j z^{-\frac{n}{2}} f(z)
$$
also has all its zeros on the unit circle.  The zeros of $f^{j+1}$
interlace the zeros of~$f^{j}$.
\end{lemma}

\begin{proof}
From 
\begin{equation}
f(z)=\sum_{0\le k\le n} a_k z^k =a_n \prod_{0\le k\le n} (z-z_j)
\end{equation}
and the fact that $1/z_j=\overline{z_j}$,
we have $f(z)=(a_0/a_n) z^n \overline{f\({1}/{\overline{z}} \)}$.
It follows that $Z(z):= (a_0/a_n)^{-\frac{n}{2}} z^{-\frac{n}{2}} f(z) $
is real on the unit circle.  By Rolle's theorem, 
$\frac{d}{d\theta} Z(e^{i\theta})$ has $n$ zeros in
$0\le\theta<2\pi$.  But these give the $n$ zeros of $Z'$,
which are the same as the zeros of~$f^1$.
\end{proof}

We now derive the analogue of Theorem~\ref{alphak}.  
Suppose the zeros of $f$ are close to equally spaced on the unit
circle:
$$
f(z)=\prod_j \(z-e\(\frac jn + \ve_j\)\)
$$
where $|\ve_j|\le\ve$ and we set $e(x)=e^{2\pi i x}$.
In the above formula and in Proposition~\ref{circlederiv}, 
unrestricted sums and products over $j$ should
be interpreted as over $j$ modulo $n$, and we consider $\ve_j$ to only
depend on $j$ modulo $n$.

Let
$$g(z)=z^{-\frac n2} f(z)
$$
and suppose $g'$ has a zero at $e(\frac kn + \frac 1{2n} + \alpha_k)$.

\begin{proposition} \label{circlederiv} In the notation above,
\begin{equation*}
\alpha_k = \frac 1{n^2} \sum_j
\frac{\ve_{j+k}}{\sin^2(\frac{\pi}{2}\,\frac{2j-1}{n})} +
O(\ve^2).
\end{equation*}

\end{proposition}

Note that the above sum is
$$
\sim \frac{4}{\pi^2} \sum_{j=-\infty}^\infty
\frac{\ve_{j+k}}{(2j-1)^2}
$$
as $n\to\infty$, so this formula matches that in
Theorem~\ref{alphak}.

\begin{corollary} $\zeta(2)=\frac{\pi^2}{6}$.
\end{corollary}

\begin{proof}
Using $\sin(x)^{-2}=x^{-2}+O(1)$, from
Proposition~\ref{circlederiv} we have
\begin{eqnarray}\label{akzetaproof}
\alpha_k &=&  \frac{4}{\pi^2} \sum_{-\frac{n}{2}<j\le \frac{n}{2}}
\frac{\ve_{j+k}}{(2j-1)^2} +O\( \frac{1}{n^2}
\sum_{-\frac{n}{2}<j\le \frac{n}{2}} |\ve_{j+k}| \)
+ O(\ve^2)\cr
&=& \frac{4}{\pi^2} \sum_{-\frac{n}{2}<j\le \frac{n}{2}}
\frac{\ve_{j+k}}{(2j-1)^2} +O\(\frac{\ve}{n}\)+O(\ve^2),
\end{eqnarray}
as $n\to\infty$ and $\ve\to 0$.

If $\ve_j=\ve$ for all $j$, then by a change of variables we see
that the zeros of $g'$ are at
$e(\frac{k}{n}+\frac{1}{2n}+\ve)$.  This means that $\alpha_k=\ve$, so
canceling $\ve$ in~(\ref{akzetaproof}) gives
$$
1=\frac{4}{\pi^2} \sum_{-\frac{n}{2}<j\le \frac{n}{2}}
\frac{1}{(2j-1)^2} +O\(\frac{1}{n}\)+O(\ve).
$$
Since $n$ and $\ve$ are arbitrary, $1=
\frac{4}{\pi^2}\sum_{j} \frac{1}{(2j-1)^2}$, which
is equivalent to $\sum_{j=1}^\infty
\frac{1}{j^2}=\frac{\pi^2}{6}$.
\end{proof}

\begin{proof}[Proof of Proposition \ref{circlederiv}]
We have
$$
\frac{g'}{g}(z)=-\frac{n}{2z}+\sum_j \frac1{z-e(\frac jn + \ve_j)},
$$
so setting $z=e(\frac kn + \frac 1{2n} + \alpha_k)$ gives
$$
\frac{n}{2 e(\frac kn + \frac 1{2n})e(\alpha_k)} = \sum_j
\frac{1}{e(\frac kn + \frac 1{2n})e(\alpha_k) - e(\frac jn )
e(\ve_j)}.
$$
Write $E(t)=e(t)-1$, so that it is easier to keep track of the
main terms:
\begin{eqnarray}
&&\frac{n}{2} e(-\frac kn - \frac 1{2n})
+\frac{n}{2} e(-\frac kn - \frac 1{2n})E(-\alpha_k)\cr
&&\phantom{XXXXXXXXXXX} = \sum_j \frac1{e(\frac kn + \frac 1{2n})
-e(\frac jn )
+ e(\frac kn + \frac 1{2n})E(\alpha_k)-e(\frac jn )E(\ve_j)} \cr
&&\phantom{XXXXXXXXXXX} = \sum_j \frac{1}{e(\frac kn + \frac
1{2n}) - e(\frac jn )} + \sum_j \frac{e(\frac jn ) E(\ve_j) -
e(\frac kn + \frac 1{2n}) E(\alpha_k)}{
        (e(\frac kn + \frac 1{2n}) - e(\frac jn ))^2} \\
        &&\phantom{XXXXXXXXXXXXX}  +O\(\sum_j\frac{|E(\ve_j)|^2+|E(\alpha_k)|^2}{
        |e(\frac kn+\frac{1}{2n})-e(\frac jn)|^3}\)\nonumber
        .
\end{eqnarray}
The first terms on both sides of the above can be canceled,
because both equal 
$$
\frac{d}{dz}\log(z^n-1)=
\frac{nz^{n-1}}{z^n-1} = \sum_{j=1}^n\frac{1}{z-e(j/n)}
$$ 
evaluated at $e(\frac kn + \frac 1{2n})$.

Rearranging to put the terms with $\alpha_k$ on one side and
multiplying by $e(\frac{ k}{n}+\frac{1}{2n})$, we have
$$
\frac{n}{2} E(-\alpha_k) +
E(\alpha_k)\sum_j\frac{e(\frac{ k}{n}+\frac{1}{2n})^2}
{(e(\frac k n +\frac{1}{2n}) - e(\frac{j}{n}))^2} 
= \sum_j \frac{e(\frac kn + \frac jn + \frac
1{2n} ) E(\ve_j)}{
        (e(\frac kn + \frac 1{2n}) - e(\frac jn ))^2} +O\(n^3\ve^2\).
$$
The sum on the left side equals $-\frac{n(n-2)}{4}$ because it equals
$$
-z^2 \frac{d^2}{dz^2}\log(z^n-1)=
\frac{nz^{2n} + (n^2-n)z^{n}}{(z^n-1)^2} = \sum_{j}\frac{z^2}{(z-e(\frac{j}{n}))^2}
$$ 
evaluated at $e(\frac kn + \frac 1{2n})$.  This brings us to
\begin{eqnarray*}
\frac{n}{2} E(-\alpha_k) - \frac14 (n^2-2n) E(\alpha_k) =&& 
\sum_j\frac{E(\ve_j)}{
        (e(\frac{k}{2n}+\frac{1}{4n}-\frac{j}{2n})-e(-\frac{k}{2n}-\frac{1}{4n}+\frac{j}{2n}))^2}
+O\(n^3\ve^2\)\\
=&&
-\sum_j\frac{E(\ve_j)}{4\sin^2(2\pi(\frac{k}{2n}+\frac{1}{4n}-\frac{j}{2n}))}
+O(n^3\ve^2).
\end{eqnarray*}
Now use the approximation $E(t)= t + O(t^2)$, for $t$ small, and solve for
$\alpha_k$:
\begin{eqnarray*}
\alpha_k &=& 
\frac 1{n^2}
\sum_j \frac{\ve_j}{\sin^2(2\pi(\frac{k}{2n} + \frac{1}{4n} - \frac{j}{2n}))}  + O(n\ve^2) \\
&=& \frac 1{n^2} \sum_j
\frac{\ve_{j+k}}{\sin^2(\frac{\pi}{2}\,\frac{2j-1}{n})}  +
O(n\ve^2) ,
\end{eqnarray*}
as claimed.
\end{proof}

\subsection{Rate of convergence to equal spacing}\label{rate}

Now we analyze the rate at which differentiation evens out 
zero spacings.

For both the midpoint and differentiation process, applying the
operation \emph{twice} has somewhat nicer properties than just
applying it once.  For example, if the points are equally spaced
then applying the operation twice leaves the sequence unchanged.
In the discussion below we will generally work with the
second derivative and second midpoint process, although 
the same ideas apply to the basic process.

We will view differentiation and midpoint averaging as two
examples of a general procedure.  Suppose $P$ is a probability
measure on the integers such that $P$ is even, $P(0)$ is the maximum, 
and $P(j)$ is increasing for $j<0$ and decreasing for $j>0$.
The conditions on $P$ are natural if one is thinking in terms of
using $P$ to smooth out the irregularities of a sequence.

Given a sequence $(\ve_n)$ define a sequence of sequences by  $\beta^0_n=\ve_n$ and
\begin{equation}
\beta^{\ell+1}_k = \sum_n P(n) \beta^\ell_{n+k}.
\end{equation}
For the second derivative process $P(0)=\frac13$ and 
$P(n)=\frac{2}{\pi^2 n^2}$ otherwise,
and for the (second) midpoint averaging process $P(0)=\frac12$,
$P(-1)=P(1)=\frac14$, and $P(n)=0$ otherwise.

In Theorem~\ref{oneoverj} we start with $z_k=k+\ve_k$, so in the
notation above 
the discrepancy of $z^\ell_{n+1}$ and  $z^\ell_{n}$ from the average 
spacing is $\beta^\ell_{k+1} - \beta^\ell_k$.
We wish to estimate this in terms of~$\ve_n$.
We can write $\beta^\ell_k$ in terms of $\ve_n$ as
\begin{equation}
\beta^{\ell}_k = \sum_n P^\ell(n) \ve_{n+k},
\end{equation}
say.
We can estimate the difference as
\begin{eqnarray}
\beta^\ell_{k+1} - \beta^\ell_k &=& \sum_n \ve_{n+k} (P^\ell(n) - P^\ell(n-1)) \nonumber \\
&\le& \ve \sum_n |P^\ell(n) - P^\ell(n-1)| \nonumber\\
&=& 2 \ve P^\ell(0) .
\label{pineq}
\end{eqnarray}
The last step requires that  $P^\ell(0)$ is the maximum, and
$P^\ell(n)$ is increasing for $n<0$
and decreasing for $n>0$.  We will see that this is the case.

It can be seen that
$P^\ell=P*P*\cdots*P$ is just the $\ell$-fold iterated
convolution of $P$ with itself, where
\begin{equation}
P*Q(n)=\sum_m P(m)Q(n-m) .
\end{equation}
The $\ell=2$ case was demonstrated in Section~\ref{almostline},
where $\beta_k$ was calculated from~$\alpha_j$.
It is straightforward to check that if $0$ is the global maximum
for $P$ and $Q$,  and both functions are  
increasing for $n<0$ and decreasing for $n>0$,
then the same holds for $P*Q$.
Thus, inequality (\ref{pineq}) applies and
we need only evaluate $P^\ell(0)=P*P*\cdots*P(0)$.

The evaluation can be done using  the properties  of convolutions and
Fourier transforms.  Given a sequence $S$, let
\begin{equation}
F_S(x)=\sum_n S(n)e^{2\pi i n x}.
\end{equation}
Then
\begin{eqnarray}
F_S(x) F_T(x)&=&\(\sum_n S(n)e^{2\pi i n x}\)\(\sum_j T(j)e^{2\pi i j x}\) \cr
&=& \sum_n \sum_j S(n)  T(j) e^{2\pi i( n+j) x}   \cr
&=& \sum_m e^{2\pi i m x} \sum_n S(n) T(m-n) \cr
&=&  \sum_m (S*T) (m)  e^{2\pi i m x} .
\end{eqnarray}
That is, convolution of sequences corresponds to multiplication of Fourier series.

The final ingredient is to note that if
$$
F(x)=\sum _j c_j e^{2\pi i j x} ,
$$
then
$$
c_0 = \int_0^1 F(x) \, dx ,
$$
provided we can integrate term-by-term,
because $\int_0^1 e^{2\pi i j x}\, dx = 0$ unless $j=0$.
Thus,
\begin{eqnarray}
 P^\ell(0)&=&P*P*\cdots*P(0) \cr
&=&  \int_0^1 \(\sum_n P(n)e^{2\pi i n x}\)^\ell dx.
\end{eqnarray}
Now we need only identify the function in the integrand above.

In the case of the (second) midpoint process, we  have
$$
\sum_n P(n)e^{2\pi i n x} = \frac12 + \frac12 \cos(2\pi x).
$$
You can either recognize the Beta integral or ask
a computer algebra package to verify that
\begin{eqnarray}
 P^\ell(0) &=& \int_0^1 \(\frac12 + \frac12 \cos(2\pi x)\)^\ell  dx\cr
&=& \frac{2\Gamma(\frac{3}{2} + \ell)}
{\left( 1 + 2\ell \right) {\sqrt{\pi }}\Gamma(1 + \ell)}  \cr
&\sim& \frac{1}{\sqrt{\pi \ell}} \\
&=&O(\ell^{-\frac12}), \nonumber
\end{eqnarray}
which is the error term from Theorem~\ref{sqrtj}, assuming~$\ve_n\ll 1$.

In the case of the (second) differentiation process, 
it follows from the Fourier expansion
$$
4\(x-\frac12\)^2=\frac13 + \frac{4}{\pi^2}\sum_{n=1}^\infty 
\frac{\cos(2\pi n x)}{n^2},
\ \ \ \ \ \ \ \ \ \ \ 
0\le x\le 1 ,
$$
that
$$
\sum_n P(n)e^{2\pi i n x} = 4 \(x-\frac12\)^2,
$$
for $0\le x\le 1$.  Thus
\begin{eqnarray}
 P^\ell(0) &=& \int_0^1 \(4 \(x-\frac12\)^2\)^\ell  dx\cr
&=& \frac{1}{1+2\ell}\\
&=& O(\ell^{-1}),  \nonumber
\end{eqnarray}
as claimed in Theorem~\ref{oneoverj}.  

The discrepancy from equal spacing is much smaller for the differentiation
process than for the midpoint process, and this is not surprising.
Both processes average among neighboring points, and that averaging
will be more effective if takes place over a larger range.  
In the setup we have described, if $P$ has finite variance 
$\sigma^2$
then $F_P(x)\sim 1- 2\pi^2\sigma^2  x^2$ for $x$ near $0$.  Therefore
$$
\int_{-\delta}^\delta F_P(x)^\ell \, dx \sim  \frac{1}{\sqrt{ 2\pi\, \sigma^2\,\ell}} ,
$$
as we saw for the midpoint process.
For the differentiation process we found that $F_P(x)=1-b |x|+O(x^2)$
for $x$ near $0$, where $b=4$.
Therefore
$$
\int_{-\delta}^\delta F_P(x)^\ell \, dx \sim  \frac{2}{b\,\ell} ,
$$ 
as in Theorem~\ref{oneoverj}.  Many other behaviors are possible.
If
$P$ has infinite variance but
decreases like  $n^{-A}$
then
$F_P(x)\sim 1-a |x|^{\alpha}$  as $x\to0$, where $\alpha=A-1$,
and so
$$
\int_{-\delta}^\delta F_P(x)^\ell \, dx \sim  \frac{2 \Gamma(\frac{1}{\alpha})}{\alpha (a\,\ell)^{\frac{1}{\alpha}}} .
$$ 
Thus, we have averaging processes where the $\ell$th iterate approaches
equal spacing as fast as any given power of~$\ell$.

It should be noted that if one violates the conditions that  $P(0)$ is a maximum and 
$P(n)$ increases (decreases) for $n<0$ ($n>0$), then the ``averaging process''
may not lead to equal spacing.  For example, if $P(-1)=P(1)=\frac12$ then 
the process leads to two interlaced sequences which separately approach
equal spacing.

\section{Proof of Theorem~\ref{cosinetheorem}}\label{cosineproof}

We wish to show that if $f$ is an entire function that meets suitable
additional conditions, then the $n$th derivative $f^{(n)}$, appropriately
rescaled, approaches the cosine function.  The real issue here is
proving that the zeros approach equal spacing with sufficient uniformity.
For if the zeros approach equal spacing, then the Hadamard product for $f^{(n)}$
can be seen to be close to $Ae^{Bx}\cos(C x+D)$ for some $A,B,C,D$.
So we first show that the zeros of $f^{(n)}$ approach equal spacing,
and then the proof is almost immediate.

Our method makes extensive use of the midpoint process we studied in 
this paper. It would be interesting to find a more direct proof,
which could possibly lead to a stronger result.
At the end of this section we also discuss the case that the
zeros of $f$ lie near, but not necessarily on, a line.

\subsection{Repeated differentiation leads to equal spaced zeros}

We have shown that under repeated differentiation the small gaps 
between zeros are becoming larger, and the large gaps are becoming smaller,
but it does not trivially follow that those gaps are  approximately
equal.  Since we know that the midpoint process gives gaps which approach
equal spacing, one possibility is to show that differentiation is
better  than midpoint at evening out the sequence:

\begin{conjecture} \label{diffisbetter}{\bf Differentiation is better than midpoint. }
Suppose $f$ is an entire function of order~1 which is real on the real
axis and has only real zeros, and suppose $z_j$ are the zeros of~$f$, listed
in increasing order.
If $q<p$ are consecutive zeros of $f'$ then 
$$
\inf \frac12 (z_{n+2}-z_n)\le p-q \le \sup \frac12 (z_{n+2}-z_n).
$$
\end{conjecture}

A similar result should also hold when the spacing between zeros of
$f$ varies slowly, in analogy to the relation between Theorem~\ref{thm:order1version} and
Theorem~\ref{thm:localversion}.  
Note that we interpret the right side of the above
inequality as ``$\infty$'' if $f$ has only finitely many zeros.

In contrast to Theorems~\ref{thm:order1version} and~\ref{thm:localversion}, the conjecture
is not true if instead we assume $p$ and $q$ are zeros of
$f'+a f$,  for if $a$ is large then $p$ and $q$ are very close to
zeros of~$f$.  Also, the conjecture is not true if $ \frac12 (z_{n+2}-z_n)$
is replaced by $ \frac13 (z_{n+3}-z_n)$.

We now show that repeated differentiation leads to equal spaced zeros.  
The proof makes use of the fact that the 
midpoint process leads to equal
spacing, but our approach is somewhat less elegant than would follow from the
above conjecture.

\begin{proposition}  Suppose $f$ is an entire function which is real on the
real line, has only real zeros, and $n_+(r)\sim n_-(r)\sim r$.  Then if
$(z^{j}_n)$ are the ordered zeros of the $j$th derivative $f^{(j)}$ then
$z^{j}_{n+1} - z^{j}_n=1+o(1)$ as $j\to\infty$.  In addition, 
$|z^{j}_n|\gg n$.
\end{proposition}

\begin{proof}
The same method as in the proof of 
Theorems~\ref{thm:order1version} and ~\ref{thm:localversion} 
shows that if 
$(z_n)$  are the ordered zeros of $f$ and
$(z'_n)$ are the ordered zeros of $f'$,
then 
\begin{equation}\label{Nspacing}
 \inf (z_{n+N}-z_n) \le z'_{n+N}-z'_n \le \sup (z_{n+N}-z_n),
\end{equation}
for any~$N$, and similarly when considering zeros with $|z_n-z_n'|<X$.
When $N=2$, the difference $z_{n+2}-z_n$ is twice the gap between
the midpoints of
consecutive zeros of~$f$.  By Theorem~\ref{sqrtj}, iterating the 
midpoint process gives sequences which approach equal spacing.
Thus, the sequence of next nearest neighbors of zeros of
$f^{(j)}$ approaches equal spacing.  That is, the zeros
of $f^{(j)}$ consist of two interlaced sequences which separately
are approaching the (same) equal spacing.

It remains to show that the whole sequence of  $f^{(j)}$ zeros
is approaching equal spacing. One possibility is to 
now prove that  $f^{(j)}(z)$ is approximately of the form
$Ae^{Bz}\(\cos(Cz+D)+E\)$, so then the zeros of derivatives of
 $f^{(j)}$ will approach equal spacing, as discussed in the 
paragraph after Theorem~\ref{cosinetheorem}.  But our goal is to
prove that the zeros become equally spaced and then use that as
a tool to prove the functional form.

We will show that the zeros of $f^{(j)}$ also consist of
\emph{three} interlaced sequences, each separately approaching
equal spacing.  From this it immediately follows that the 
entire sequence is approaching equal spacing.  
We will use~(\ref{Nspacing}) with $N=3$.
Note that
\begin{eqnarray}
z_{n+3}-z_n &=& (z_{n+3}+z_{n+2}+z_{n+1}) - (z_{n+2}+z_{n+1}+z_{n}) \cr
&=& 3 ({\tilde{z}}_{n+2} - {\tilde{z}}_{n+1} ) ,
\end{eqnarray}
where
$ {\tilde{z}}_{n} = \frac13 (z_{n-1}+z_n+z_{n+1}) $.
By either imitating the proof of Theorem~\ref{sqrtj} or using the
method of Section~\ref{rate} with $P(-1)=P(0)=P(1)=\frac13$, we see
that iterating the averaging process $(z_n)\to({\tilde{z}}_n)$ 
leads the sequence to approach equal spacing. 
So $z^j_{n+3}-z^j_n$ is approaching equal spacing, as required.

The final assertion follows from the fact that $z^0_n\sim n$ 
and that each successive zero set interlaces the previous one.
\end{proof}

\subsection{Proof of Theorem~\ref{cosinetheorem}}
Assume $\kappa=1$, so the zeros of $f^{(j)}$ are 1 apart on average.
We first shift the function slightly so that the zeros are arranged
conveniently.
Suppose the smallest non-negative zero of  $f^{(j)}$ is at $w_1$ and the
largest negative zero is at $w_{-1}$.  
Choose $d_j=\frac12(w_{-1}+w_1)$ and let $z_1$, $z_2$,\ldots denote the
positive zeros of $f^{(j)}(z+d_j)$ and $z_{-1}$, $z_{-2}$,\ldots the negative
zeros.   

We have the Hadamard factorization 
\begin{eqnarray}
f^{(j)}(z+d_j)&=&A_j \exp\({B_j z}\) 
\prod_{n}\(1-\frac{z}{z_n}\)\exp\(\frac{z}{z_n}\) \cr
&=& A_j \exp\({B_j' z}\)  \prod_{1\le n\le Y}\(1-\frac{z}{z_{-n}}\) \(1-\frac{z}{z_n}\)\cr
&&  \ \ \ \ \  \ \ \ \ \ \ \ \  \ \ \ \ \ \ \ \ \ \ 
\times
\prod_{|n|>Y}\(1-\frac{z}{z_n}\)\exp\(\frac{z}{z_n}\)\cr
&=& A_j \exp\({B_j' z}\) F_j(z) G_j(z),
\end{eqnarray}
say, where $Y$ is to be chosen later.
We must show that if $|z|<X$ then $F_j(z)\to \cos(\pi z)$ and 
$ G_j(z)\to 1$. 

Given $X$, $\ve>0$, first choose $Y$ so that $|G_j(z)-1| < \ve e^{-\pi X}$ for
$|z|<X$ and all $j>0$.  This is possible because $z_n\gg n$.

Let $F$ and $G$ denote the above products in the special case that
$z_{\pm n}=(n\mp \frac12)^{-1}$ for all $n>0$.  Since $F(z)G(z)=\cos(\pi z)$,
it follows that $|\cos(\pi z)-F(z)|<\ve$ for $|z|<X$.
 
For each $n$, as $j\to\infty$ the gap $z^j_n-z^{j}_{n-1}$ is approaching $1$, so 
we can a choose $j$ so that 
$z_{\pm n}$ is sufficiently close to $(n\mp \frac12)^{-1}$ for $0<n\le Y$
to ensure
that $|F_j(z)-F(z)|<\ve$ for $|z|<X$.  
Using $|\cos(\pi z)|< e^{\pi X}$ for $|z|<X$ and combining all estimates we have 
$|A_j^{-1} \exp(-B_j' z) f^{(j)}(z+D_j) - \cos(\pi z)|< 4\ve$
for $|z|<X$,
which completes the proof.

\subsection{Zeros not on a line}\label{sec:notonline}

We have seen that repeated differentiation leads to a function whose
zeros approach equal spacing.  Surprisingly, this does not always require 
the zeros of the original function to lie on a line (or circle), but only that
the zeros lie in a suitable neighborhood of the line (or circle).

This is easy to see in the case of the circle.  
Suppose $f$ is a degree $n$ polynomial 
and let $g(z)=z^{-n/2}f(z)$, so
\begin{equation}
g(z)=a_n z^{\frac{n}{2}} + a_{n-1}z^{\frac{n}{2}-1} + \cdots 
        a_0 z^{-\frac{n}{2}} .
\end{equation}
Then
\begin{equation}
z g'(z)=\frac{n}{2} a_n z^{\frac{n}{2}} + 
\(\frac{n}{2}-1\)a_{n-1}z^{\frac{n}{2}-1} + \cdots 
      \(-\frac{n}{2}\)  a_0 z^{-\frac{n}{2}} ,
\end{equation}
and the ``$k$th derivative'' equals
\begin{equation}
\(z\frac{d}{dz}\)^k g(z)=
\(\frac{n}{2}\)^k a_n z^{\frac{n}{2}} + 
\(\frac{n}{2}-1\)^k a_{n-1}z^{\frac{n}{2}-1} + \cdots 
      \(-\frac{n}{2}\)^k  a_0 z^{-\frac{n}{2}} .
\end{equation}
If $a_0\not=0$ then the factors $(\frac{n}{2})^k$ from the first and last terms dominate
everything else, so the zeros are approaching the zeros of
$z^n + (-1)^k \frac{a_0}{a_n} $.
Note that we didn't actually require the zeros of the original
polynomial $f$ to be on the unit circle.
As long as $0$ is not a root of $f$, the above process gives a sequence
of functions whose zeros eventually lie on a circle, and in fact
approach equal spacing on a circle!

For zeros on a line, we have the following result of Young-One Kim 
(this is a slight revision of Theorem~2 from~\cite{kim}):

\begin{theorem}\label{kimtheorem}
Let $f(z)$ be a nonconstant real entire function, $0<\rho\leq 2$, 
and assume that $f(z)$
is of order less than $\rho$ or is of order $\rho$ and minimal type. If there is a
positive real number $A$ such that all the zeros of $f(z)$ are distributed
in the infinite strip $|\Im z|\leq A$, then for any positive constant $B$
there is a positive integer $n_1$ such that $f^{(n)}(z)$ has only real
zeros in $|\Re z|\leq Bn^{\frac{1}{\rho}}$ for all $n\geq n_1$.  
\end{theorem}

In other words, if the zeros don't start out too far from the real axis, 
then they end up on the real axis as you differentiate.  
Thus, Theorem~\ref{cosinetheorem} can be modified
to only assume that the zeros lie in a strip around the real axis,
and Conjecture~\ref{Xiconjecture} should
only require the zeros to lie near the real axis.

The theorem applies to the Riemann $\Xi$-function with
$\rho=1+\varepsilon$. 
It is interesting to note that Conrey~\cite{co} has shown that
$\Xi^{(j)}$ has a positive proportion of its zeros on the real 
axis, and that proportion is $1+O(j^{-2})$  as  $j\to\infty$.  

If $f(z)=\sum_{j=0}^\infty c_j x^j/j!$ is entire of order~1 
and all $c_j$ are real, 
then a necessary 
condition for the zeros of $f$ to be real is  the Tur\'{a}n inequalities 
$c\sb{k}\sp{2}- c\sb{k-1}c\sb{k+1}\ge 0$, for $k\ge 1$. 
The above ideas suggest that the Tur\'{a}n inequalities
should hold for sufficiently large~$k$, provided only that the zeros
of $f$ lie in a neighborhood of the real axis.

\end{document}